\RequirePackage{fix-cm} 
\documentclass[a4paper,twoside,11pt,reqno]{amsart}

\usepackage{amsaddr,xcolor}
\usepackage{indentfirst,a4wide,verbatim}
\usepackage{amsmath,amsfonts,amssymb,amsgen,amsbsy,dsfont,eucal,mathrsfs,dsfont,stmaryrd,amsthm}
\usepackage{eucal,esint,dsfont,mathrsfs,MnSymbol}
\usepackage{stmaryrd} 
\usepackage{nomencl}
\usepackage{graphicx}
\usepackage{graphics}
\usepackage{float} 
\usepackage[section]{placeins} 
\usepackage{algorithm} 
\usepackage[noend]{algpseudocode} 
\usepackage{xcolor}
\usepackage{multirow}

\usepackage{subfig}
\usepackage{dirtytalk}
\usepackage{appendix}
\usepackage{color}
\usepackage{cleveref}
\usepackage{ulem}
\usepackage{graphicx}
\usepackage[square,numbers]{natbib}
\newcommand{\ud}{\mathrm{d}}

\newcommand{\ue}{\mathrm{e}}
\newcommand{\ui}{\mathrm{i}}

\newcommand{\half}{{\textstyle{\frac{1}{2}}}}
\newcommand{\ihalf}{{\textstyle{\frac{\ui}{2}}}}

\newcommand{\quat}{{\textstyle{1\over4}}}

\renewcommand{\Re}{\operatorname{Re}}

\newcommand{\be}{\begin{equation}}
\newcommand{\ee}{\end{equation}}
\newcommand{\eqdef}{\stackrel{\text{\tiny{def}}}{=}} 

\title{Optimized integrating factor technique for  
	Schr{\"o}dinger-like equations}
\author{M.~Lovisetto, D.~Clamond and B.~Marcos}

\address{Universit\'e C\^ote d'Azur, CNRS UMR 7351, Laboratoire J. A. Dieudonn\'e,\\
Parc Valrose, F-06108 Nice cedex 2, France.}
\email{martino.lovisetto@univ-cotedazur.fr}

\begin{document}

\begin{abstract}
The integrating factor technique is  widely used to solve numerically 
(in particular) the Schr{\"o}dinger equation in the context of spectral 
methods. Here, we present an improvement of this method exploiting the 
freedom provided by the gauge condition of the potential. Optimal gauge 
conditions are derived considering the equation and the temporal numerical 
resolution with an adaptive embedded scheme of arbitrary order. 
We illustrate this approach with the nonlinear Schr{\"o}dinger (NLS) and 
with the Schr{\"o}dinger--Newton (SN) equations. We show that this  
optimization increases significantly the overall computational speed, 
sometimes by a factor five or more. This gain is crucial for long time 
simulations. 	
\end{abstract}

\maketitle

\medskip


\section{Introduction}\label{intro}

The Schr\"odinger equation is used in many field of Physics.
In dimensionless units, it takes the form
\begin{equation} \label{def-schro}
\ui\,\partial_t\,\psi\ +\ \half\,\nabla^2 \psi\ -\ V\,\psi\ =\ 0,
\end{equation}
where $\psi$ is a function of the spatial coordinates $\boldsymbol{r}$ and 
of the time $t$, $\nabla^2$ is the Laplace operator, and the potential $V$ 
is generally a function of space and time and, possibly, a functional of $\psi$. 
In this paper, we focus more specifically on the nonlinear Schr\"odinger 
equation ($V$ proportional to $|\psi|^2$) and on the Schr\"odinger--Newton 
(or Schr\"odinger--Poisson) equation in which the Laplacian of the potential 
is proportional to $|\psi|^2$. These special cases were chosen for clarity 
and because they are of practical interest, but the method presented here 
can be extended to more general potentials (and equations).

With the nonlinear Schr\"odinger  equation (NLS) considered here, 
the (local nonlinear) potential is  
\begin{equation}\label{potSN}
V\ =\ g\,|\psi|^2,
\end{equation}
where $g$ is a coupling constant. For $g>0$ the interaction is repulsive, 
while, for $g<0$ it is attractive.  The NLS describes various physical phenomena, such as 
Bose--Einstein condensates \citep{dalfovo1999theory}, laser beams in some nonlinear 
media \citep{hasegawa1989optical}, water wave packets \citep{kharif2008rogue}, etc.

With the Schr\"odinger--Poisson equation, the potential is given by the Poisson 
equation 
\begin{equation}\label{potNLS}
\nabla^2\,V\ =\ g\, |\psi|^2,
\end{equation}
where $g$ is another coupling constant, the interaction being attractive if $g>0$ 
and repulsive if $g<0$. It is therefore {\it nonlinear\/} and {\it non-local}, 
giving rise to collective phenomena \citep{binney}, appearing for instance in optics  
\citep{dabby1968thermal,rotschild2005solitons,rotschild2006long}, Bose--Einstein 
condensates  \citep{guzman2006gravitational}, cosmology 
\citep{hu2000fuzzy,paredes2016interference,marsh2019strong} and theories describing 
the quantum collapse of the wave function \citep{diosi2014gravitation,penrose1996gravity}. 
It is also used as a model to perform cosmological simulations 
in the semi-classical limit \citep{widrow1993using}. 

The above equations cannot be solved analytically (except for very special cases) and 
numerical methods must be employed. In this paper, we focus on spectral methods for the 
spatial resolution, i.e., methods that are based on fast Fourier transform (FFT) techniques, 
that are specially efficient and accurate \citep{CanutoEtAl2006-1}. For the temporal 
resolution, two families of methods are commonly employed to solve Schr\"odinger-like 
equations: integrating factors \citep{lawson1967generalized} 
and split-step integrators \citep{blanes2008splitting}. 
The latter methods have been used to integrate both the SN and NLS equations, but the 
former is used essentially to solve the NLS, with very performing results 
\citep{bader2019efficient,besse2017high}.
In this note, we focus on the former technique, which consists in integrating analytically 
the linear part of the equation and integrating numerically the remaining nonlinear part 
with a classical method \cite{lawson1967generalized}. The principle of the method is described as 
follow. 

Writing the Schr\"odinger equation in the generic form
\begin{equation}
\ui\,\partial_t\,\psi\ =\ F\!\left(\boldsymbol{r},t,\psi\right), \qquad
\psi\ =\ \psi\!\left(\boldsymbol{r},t\right),
\end{equation}
the right-hand side is split into linear and nonlinear parts 
\begin{equation}    \label{eqIF1}
\ui\,\partial_t\,\psi\ +\ \mathcal{L}\,\psi\ =\ \mathcal{N}\!\left(\boldsymbol{r},t,\psi\right),
\end{equation}
where $\mathcal{L}$ is an easily computable autonomous linear  
operator and $\mathcal{N}\eqdef F+\mathcal{L}\psi$ is the remaining 
(usually) nonlinear part.
At the $n$-th time-step, with $t\in[t_{n},t_{n+1}]$, considering the change of 
dependent variable
\begin{equation}
\phi\ \eqdef\ \exp\!\left[\,(t-t_n)\,\mathcal{L}\,\right]\psi \qquad\implies\qquad 
\ui\,\partial_t\,\psi\ =\ \exp\!\left[\,(t_n-t)\,\mathcal{L}\,\right]\left(\,\ui\,\partial_t\,\phi\,-\,\mathcal{L}\,
\phi\,\right),
\end{equation}
so $\phi=\psi$ at $t=t_n$, the equation \eqref{eqIF1} is rewritten
\begin{equation} \label{eqIF2}
\ui\,\partial_t\,\phi\ =\,\exp\!\left[\,(t-t_n)\,\mathcal{L}\,\right]\mathcal{N}.
\end{equation}
The operator $\mathcal{L}$ being well chosen, the stiffness of \eqref{eqIF1} is considerably 
reduced and the equation \eqref{eqIF2} is (hopefully) well approximated by algebraic polynomials 
for $t\in[t_n;t_{n+1}]$. Thus, standard time-stepping methods, such as an adaptive 
Runge--Kutta method \cite{alexander1990solving,butcher2016numerical}, can be used to efficiently 
solve \eqref{eqIF2}. 
To do so, the solution is evaluated at two different orders and a local error is estimated 
as the difference between those quantities. Popular integrators can be found in 
\citep{alexander1990solving,dormand1980family}.

It is straightforward to apply this strategy to the Schr{\"o}dinger equation \eqref{def-schro} 
since  $\half \nabla^{2}\psi$ and the potential $V$ are, respectively, 
linear and nonlinear operators of $\psi$. By switching to Fourier space in position, the 
equation becomes
\begin{equation}\label{schro-fourier}
\ui\,\partial_t\,\widehat{\psi}\ -\ \half\,k^2\,\widehat{\psi}\ -\ \widehat{\/V\,\psi\/}\ =\ 0,
\end{equation}
where \say{hats} denote the Fourier transform of the underneath quantity and $k\eqdef|\boldsymbol{k}|$ 
($\boldsymbol{k}$ the wave vector). The equation is now in a form where the application of the 
integrating factor (IF) method is straightforward, i.e., \eqref{schro-fourier} becomes
\begin{equation}\label{SNIF1}
\ui\,\partial_t\,\phi\ =\ -\/\ui\,\exp[\ihalf\/k^2\/(t-t_n)\/]\,\widehat{\/V\,\psi\/},
\end{equation}
where $\phi(\boldsymbol{k},t)\eqdef\widehat{\psi}(\boldsymbol{k},t)\exp[\ihalf\/k^2\/(t-t_n)\/]$. 
If the nonlinear part of the equation is zero, then $\ui\,\partial_t\,\phi=0$ and any (reasonable) temporal 
scheme will produce the exact solution $\phi(t)=\phi(t_n)$. In other words, the integrating 
factor technique is exact for linear equations. This indicates that the numerical errors depend 
on the magnitude of the nonlinear part. Therefore, in order to minimise these errors, a strategy 
consists in minimizing the magnitude of $\mathcal{N}$ at each time-step.
To do so, we exploit the gauge invariance of the Schr\"odinger equation: if $\psi$ is a solution of \eqref{def-schro} at a given time $t$, then $\Psi\ \eqdef \psi\,\ue^{-\/\ui\,\mathcal{C}\,t}$ is a solution of
\begin{equation} \label{eq-gauge-rotation}
\ui\,\partial_t\,\Psi\ +\ \half\,\nabla^2 \Psi\ -\ \left(V + \mathcal{C} \right)\,\Psi\ =\ 0,
\end{equation}
as one can easily verify. Thus, at each time-step, adding a constant $\mathcal{C}_n$ to $V$ in \eqref{def-schro}, modifies the solution as
\begin{equation} \label{gauge-rotation}
\psi(t_n)\ \to\ \psi(t_n)\,\ue^{-\/\ui\/\varphi}, \qquad
\varphi\ \eqdef\ \sum_{j=1}^{n} \mathcal{C}_j\,h_j,
\end{equation}
where $h_j\eqdef t_{j+1}-t_j$ is the $j$-th time-step. 
Of course, at the end of the computations, the operation \eqref{gauge-rotation} can be easily reverted 
if the original phase is relevant. Using this procedure, we observed up to a five-fold speed increase 
(the overall computational time is divided by about five) compared to taking $\mathcal C_n=0$. Of course, 
the speed-up varies depending on the initial condition, of the (spatial and temporal) numerical schemes and on the choice of gauge corresponding to $\mathcal C_n=0$.

In this paper, we derive some analytic formulas giving an optimal 
$\mathcal{C}_n$ in order to maximise the time-step, i.e., to minimize 
the overall computational time of the numerical resolution. We emphasize 
that the thus obtained optimal values of $\mathcal{C}_n$ do not affect 
the accuracy of the numerical solution, leaving it unchanged with respect
to the $\mathcal{C}_n=0$ case. Two strategies are presented. 
In section \ref{sect-C}, a first `natural' approach to derive a suitable 
$\mathcal{C}_n$ is based on the analytical structure of the equation and 
it is independent of the numerical algorithm employed for its resolution.
More precisely, $\mathcal{C}_n$ is obtained minimizing a norm of the 
right-hand side of the equation \eqref{eqIF1}. 
This provides an easy manner to obtain a formula that is moreover computationally 
`cheap'. This expression is however only near optimal, so a `better' 
expression is subsequently derived. 
Considering both the equations and the numerical algorithms, a second 
optimal expression for $\mathcal{C}_n$ is derived in section \ref{secoptCn}. 
This approach consists in minimizing exactly the numerical 
error and thus explicitly dependents on the numerical scheme. This provides  
a more accurate, but computationally expensive, solution. 
The advances of these special choices are illustrated numerically in 
section \ref{secnum}. Finally, a summary and perspectives are drawn in 
section \ref{sect-conclusions}.

\section{Near optimal $\mathcal C_n$} \label{sect-C}

As mentioned above, if properly chosen, the integrating factor is able to reduce the stiffness of the equation, making the numerical integration more efficient. In addition, the magnitude of the nonlinear part of \eqref{eqIF2} also contributes to the efficiency of the numerical integration. Specifically, if $\mathcal{N}$ is zero, $\partial_t\,\phi=0$ and the integrating factor technique is exact. Thus, the efficiency of the algorithm is expected to increase as the magnitude of $\mathcal{N}$ gets 
smaller, and subsequently the overall computational time should 
be reduced. Here, we show how to choose the arbitrary constant $\mathcal{C}_n$ in order to reduce the magnitude of the nonlinear part $\mathcal{N}$. 
In the case of the Schr{\"o}dinger equation, we have
\begin{align}
	\mathcal{N}\!\left(\boldsymbol{k},t;\phi;\mathcal{C}_n\right)\ =\  -\/\ui\,\exp[
	\ihalf\/k^2\/(t-t_n)\/]\,\mathcal{F}\{\left(V+\mathcal{C}_n\right)\,\psi \},
\end{align}
where $\mathcal{F}$ denotes the Fourier transform and $\psi(\boldsymbol{x},t)= 
\mathcal{F}^{-1}\{\exp[-\ihalf\/k^2\/(t-t_n)\/]\,{\phi}(\boldsymbol{k},t)\}$. 

A natural strategy is to minimize the $L_2$-norm, namely
\begin{align}
	\mathcal{G}_n(\mathcal{C}_n)\eqdef\frac{1}{M}\,\sum_{m=-\left[M/2\right]}^{\left[M/2\right]-1}|\,\mathcal{N}\!\left(\boldsymbol{k}_m,t_n;\phi;\mathcal{C}_n\right)|^2,
\end{align}
where $M$ is the number of spatial modes, square brackets denote the integer part and $\boldsymbol{k}_m$ is the $m$-{th} Fourier mode. The explicit expression of $\mathcal{G}_n$ can be found exploiting the definition of the discrete Fourier transform. For simplicity, we do the calculations 
in one dimension (1D) without loss of generality, since the final result is independent of the spatial 
dimension $d$. From Parseval theorem, one obtains 
\begin{align}
	\mathcal{G}_n(\mathcal{C}_n)\ 
	&= \sum_{\ell=-\left[M/2\right]}^{\left[M/2\right]-1} (V_\ell+\mathcal{C}_n)^2\,|\psi_{\ell}|^2,
\end{align}
where $\psi_{\ell} \eqdef \psi({\boldsymbol x}_{\ell})$ and $V_{\ell} 
\eqdef V({\boldsymbol x}_{\ell})$ at time $t_n$.
Since the function $\mathcal{G}_n(\mathcal{C}_n)$ is a second-order polynomial in $\mathcal{C}_n$, it admits an unique minimum, which is obtained from the equation $\ud G_n(\mathcal{C}_n)/\ud\mathcal{C}_n=0$, yielding
\begin{equation}\label{L2normres}
\mathcal{C}_n\,=\,-\left(\,\sum_{\ell=-\left[M/2\right]}^{\left[M/2\right]-1} V_\ell \,\left|\psi_\ell\right|^2\,\right)\,\left/\,\left(\,\sum_{\ell=-\left[M/2\right]}^{\left[M/2\right]-1} \left|\psi_\ell\right|^2\,\right)\right.\eqdef\,\tilde{\mathcal{C}}_n.
\end{equation}
Therefore, at each time step $n$, $\tilde{\mathcal{C}}_n$, which is the 
value of $\mathcal{C}_n$ minimizing the $L_2$-norm of $\mathcal{N}$, is 
obtained from \eqref{L2normres}. We show below that even though this approach 
is not unique (i.e., different norms could be considered), the provided 
solution is quite advantageous compared to 
others, being computationally cheap and independent on the order of the 
numerical scheme.

\section{Optimal $\mathcal C_n$}\label{secoptCn}

We show here another way to choose the arbitrary constant $\mathcal C_n$ 
in order to improve the algorithm efficiency and reduce the overall 
computational time. This approach is based on the principles of the adaptive 
time-step procedure, where at each time step $n$, an error $\Delta_n$ between 
two approximated solutions of different orders is estimated. Since the smaller 
this quantity the larger the time-step, minimizing $\Delta_n$ allows to choose 
a larger time-step, speeding-up the numerical integration and keeping roughly 
the same numerical error. More specifically, the error $\Delta_n$ depends on 
the arbitrary constant $\mathcal{ C}_n$, hence the minimization can be 
performed (see below) choosing an appropriate $\mathcal{C}_n$. 
In this section, we first recall the method for determining the size of the 
time step used in the Runge--Kutta procedures; interested readers should 
refer to \cite{alexander1990solving} for further details. Although the 
determination of $\mathcal{C}_n$ can be formally presented for any 
embedded Runge--Kutta schemes, this results in very cumbersome calculations 
with little insights. Thus, for brevity and clarity, we illustrate the method 
with the Heun method (that is a second-order Runge--Kutta method with an 
embedded first-order explicit Euler scheme for the time stepping 
\cite{suli2003introduction}). We then sketch-out how this procedure can be 
implemented for generic embedded Runge--Kutta methods.

\subsection{Principle of the adaptive time-step procedure}

For the time stepping, embedded Runge--Kutta methods estimate the quadrature 
error comparing the results of two orders of the time integrator 
\cite{alexander1990solving}. For a solver of order $N$ with an embedded 
$(N-1)$-order scheme (hereafter schemes of orders $\{N,N-1\}$), at the 
$n$-th time step, the error $\Delta_n$ is \cite{wanner1996solving}
\begin{equation}\label{delta_n}
\Delta_n\ \eqdef\ \sqrt{\frac{1}{M} \sum_{m=-\left[M/2\right]}^{\left[M/2\right]-1} \left( \frac{\left|\,\phi(\boldsymbol{k}_m,t_n)\,-\,
		\widetilde{\phi}(\boldsymbol{k}_m,t_n)\,\right|}{\mbox{\sc Tol}\,+\,\max\!\left(\left|\phi(\boldsymbol{k}_m,t_n)\right|,
		\left|\tilde{\phi}(\boldsymbol{k}_m,t_n)\right|\right)\times\mbox{\sc Tol}}
	\right)^2 },
\end{equation}
where $M$ is the number of spatial modes, square brackets denote the integer part, $\phi(\boldsymbol{k}_m, t_n)$ is the $N$-{th} order solution at 
the $m$-{th} Fourier mode, the \say{tilde} notation indicating the solution at order $N-1$, and $\mbox{\sc 
	Tol}$ is the tolerance (parameter defining the desired precision of the time-integration). The time step $h_n$ 
is accepted if the error $\Delta_n$ is smaller than the tolerance $\mbox{\sc Tol}$, otherwise $h_n$ 
is reduced and this step is recomputed. $h_n$ being accepted, the next time step $h_{n+1}$ is obtained 
assuming the largest error equal to the tolerance. In order to avoid an excess of rejected time steps, 
we use the Proportional Integral (PI) Step Control \cite{wanner1996solving}, which 
chooses the optimal time step $h_{n+1}$ as 
\begin{equation}\label{h_opt}
h_{n+1}\ =\ h_{n}\,\Delta_n^{\,-b}\,\Delta_{n-1}^{\,c},
\end{equation}
where $b=0.7/p$, $c=0.4/p$, $p$ being the order of the chosen integrator \cite{gustafsson1991control}. 
Interested readers should refer to \cite{alexander1990solving} for details on this classical procedure.

\subsection{Optimum time step}

Since the constant $\mathcal{C}_n$ can be chosen freely, we seek for the value of $\mathcal{C}_n$ providing 
the largest $h_{n+1}$, namely, to maximize the right-hand side of \eqref{h_opt}. 
Since $h_n$ and $\Delta_{n-1}$ are determined at the previous time-step, only $\Delta_n$ in \eqref{h_opt} depends on $\mathcal{C}_{n}$. Thus, in order to maximize  $h_{n+1}$,  
$\Delta_{n}$ must be minimized, i.e., one must solve $\ud \Delta_{n}/\ud\mathcal{C}_n
=0$. 
This derivation being characterized by cumbersome algebra for general 
embedded Runge--Kutta schemes, we illustrate the case of the Heun algorithm 
(that is a second-order Runge--Kutta method with an embedded first-order 
explicit Euler scheme for the time stepping \cite{suli2003introduction}), 
the principle being the same for higher order integrators. Also for 
simplicity, we give the calculations in one dimension (1D) without loss 
of generality, since the final result is independent of the spatial 
dimension.

\subsubsection{Optimum $\mathcal{C}_n$ for Heun's method} \label{sect-heun}	

Heun's method consists, here, in solving the initial value problem (for $t\geqslant t_n$)
\begin{align}
	\ui\,\partial_t\,\phi\ &=\ f({\boldsymbol k},t;\phi;\mathcal{C}_n)\, \eqdef\, -\/\ui\,\exp\!\left[\,\ihalf\,k^2\,(t-t_n)\,\right]\mathcal{F}\{\,(V+\mathcal{C}_n)\,\psi\,\},
\end{align}
and
\begin{align}
	\phi({\boldsymbol k},t)\, &\eqdef\,  \exp\!\left[\,\ihalf\,k^2\,(t-t_n)\,\right]\mathcal{F}\{
	\psi({\boldsymbol x},t)\}.
\end{align}
Hereafter, for brevity, we denote
\begin{align}
	\phi_n\, =\, \phi_n({\boldsymbol k})\ \eqdef\ \phi({\boldsymbol k},t_n), \qquad \psi_n\, =\, \psi_n({\boldsymbol x}) \eqdef\ \psi({\boldsymbol x},t_n), \qquad V_n\,=\,V_n({\boldsymbol x})\eqdef V({\boldsymbol x},t_n).
\end{align}
At time $t=t_{n+1}$, the first- and second-order (in $h_n$) approximations of $\phi$, 
respectively $\widetilde{\phi}_{n+1}$ and $\phi_{n+1}$, are
\begin{align}
	\widetilde{\phi}_{n+1}\ &=\ \phi_n\ +\ h_n\,f({\boldsymbol k},t_n;\phi_n;\mathcal{C}_n), \\
	\phi_{n+1}\ &=\ \phi_n\ +\ \half\,h_n\left[\,f({\boldsymbol k},t_n;\phi_n;\mathcal{C}_n)\,+\,
	f\!\left({\boldsymbol k},t_n+h_n;\phi_n+h_n\/f({\boldsymbol k},t_n;\phi_n;\mathcal{C}_n);\mathcal{C}_n\right)\,\right].
\end{align}
The next time-step $h_{n+1}$ is chosen using equation \eqref{h_opt}. For our equation, the difference 
between the first- and second-order approximations $\Delta\phi_{n+1}\eqdef\left| \phi_{n+1} 
- \widetilde{\phi}_{n+1} \right|$ is such that
\begin{align}\label{deltaphi}
	\left( \Delta \phi_{n+1}\right)^2\ &=\ \quat\,h_n^{\,2}\left|\,f({\boldsymbol k},t_n;\phi_n;\mathcal{C}_n)\,-\,
	f\!\left({\boldsymbol k},t_n+h_n;\phi_n+h_n\/f({\boldsymbol k},t_n;\phi_n;\mathcal{C}_n);\mathcal{C}_n\right)\,\right|^2\nonumber\\
	& =\  \quat\,h_n^{\,2}\left| f({\boldsymbol k},t_n;\phi_n;\mathcal{C}_n)\,+\, \ui\,\ue^{\ui\/k^2\/h_n/2}
	\times\right.\nonumber\\ 
	&\qquad\qquad \left.
	\mathcal{F}\left\{(V_{n+1}+\mathcal{C}_n)\/\mathcal{F}^{-1}\left\{\ue^{-\ui\/k^2\/h_n/2}\left(\phi_n+h_n 
	f({\boldsymbol k},t_n;\phi_n;\mathcal{C}_n)\right)\right\}\right\}\/\right|^2,
\end{align}
where $V_{n+1}=V({\boldsymbol x},t_n+h_n)$.
We note that the absolute value in \eqref{deltaphi} is of first-order 
in $h_n$, as one can easily check with a Taylor expansion around 
$h_n=0$, so $\left( \Delta \phi_{n+1}\right)^2=O\!\left(h_n^{\,4}
\right)$. More precisely, after some elementary algebra, one finds 
\begin{align}
	\left( \Delta \phi_{n+1}\right)^2\ =\  &\quat\,h_n^{\,4}\,\,\left| \mathcal{F}\left\{(V_n+
	\mathcal{C}_n)^2\,\psi_n \right\}\ +\ \ui\,\mathcal{F}\left\{\partial_t V_n\,\psi_n \right\}\ +\ \mathcal{F}\left\{\partial_xV_n\,\partial_x\psi_n \right\} \right. \\ \nonumber
	&\left.+\ \half\,\mathcal{F}\left\{\partial_{xx}V_n\,\psi_n \right\}\ \right|^2\ +\ O\!\left(h_n^{\,5}\right),
\end{align}	
which, defining
\begin{align}
	\alpha({\boldsymbol x}, t;\mathcal{C}_n)\, &\eqdef\, (V_n+
	\mathcal{C}_n)^2\,\psi_n, \qquad
	\beta({\boldsymbol x}, t)\, \eqdef\, \ui\,\partial_t V_n\,\psi_n\ +\ \partial_xV_n\,\partial_x\psi_n +\ \half\,\partial_{xx}V_n\,\psi_n, 
\end{align}
can be rewritten as
\begin{align}\label{Deltaphi4}
	\left( \Delta \phi_{n+1}\right)^2\ =\  &\quat\,h_n^{\,4}\,\,\left|
	\mathcal{F}\left\{\alpha(\mathcal{C}_n) +  \beta \right\} \right|^2\ 
	+\ O\!\left(h_n^{\,5}\right).
\end{align}	
Introducing the mean quadratic error 
\begin{equation}\label{quaderr}
E_n(\mathcal{C}_n)\,\eqdef\ 
\frac{1}{M}\,\sum_{m=-\left[M/2\right]}^{\left[M/2\right]-1} 
\Delta \phi_{n+1}^2(\boldsymbol{k}_m,t_n;\mathcal{C}_n),
\end{equation}
substituting \eqref{Deltaphi4} into \eqref{quaderr} and exploiting 
the definition of the discrete Fourier transform, one 
obtains (using Parseval theorem) 
\begin{align}\label{result-heun}
	E_n(\mathcal{C}_n)\ &=  
	\frac{1}{M}\,\sum_{m=-\left[M/2\right]}^{\left[M/2\right]-1} 
	\left|\sum_{\ell=-\left[M/2\right]}^{\left[M/2\right]-1} 
	\ue^{-\/2\/\ui\/\pi\/m\/\ell/M}
	(\alpha_{\ell}(\mathcal{C}_n)+\beta_{\ell})\,\right|^2\,
	\quat\,h_n^{\,4}\ +\ O\!\left(h_n^{\,5}\right)  \nonumber\\
	&= \sum_{\ell=-\left[M/2\right]}^{\left[M/2\right]-1} \left|\alpha_{\ell}(\mathcal{C}_n) +  \beta_{\ell}\right|^2\,\quat\,h_n^{\,4}\,+\ O\!\left(h_n^{\,5}\right).
\end{align}
The minimum of $E_n(\mathcal{C}_n)$, obtained from the equation $\ud E_n(\mathcal{C}_n)/\ud\mathcal{C}_n=0$, is such that
\begin{align}\label{heunresult}
	\sum_{\ell=-\left[M/2\right]}^{\left[M/2\right]-1} \frac{\ud \left|\alpha_{\ell}(\mathcal{C}_n)\right|^2}{\ud \mathcal{C}_n}\ +\  2\,\Re \left( \frac{\ud  \alpha_{\ell}(\mathcal{C}_n)}{\ud \mathcal{C}_n} \beta^*_{\ell} \right) \,=\,0,
\end{align}
Therefore, the optimum $\hat{\mathcal{C}}_n$ providing 
the largest $h_{n+1}$, in the case of Heun's method, is a  solution of \eqref{heunresult}.

\subsubsection{Optimum $\mathcal{C}_n$ for generic embedded 
	Runge--Kutta schemes}\label{sect-general}

The optimum $\mathcal C_n$ for general embedded Runge--Kutta schemes can be 
obtained following the same principles illustrated above with the Heun 
algorithm. However, the algebraic calculations get rapidly very cumbersome, 
leading to expensive computations that, in most cases, exceeds the time gained 
with a larger step. Here, we sketch-out the procedure for generic embedded 
Runge--Kutta methods, considering solvers of order $N$ with an embedded 
$(N-1)$-order scheme (for other embedded or extrapolation methods, the 
procedure is completely analogue). 
For a $s$-stage method, the error $\Delta \phi_{n+1}$ can be written as \citep{wanner1996solving}
\begin{align}
	\left( \Delta \phi_{n+1}\right)^2\,=\,\left| \sum_{\ell=1}^s d_{\ell}\, w_{\ell}\right|^2, \quad d_{\ell} \eqdef a_{s,\ell}-b_{\ell}, \quad 	w_{\ell} \eqdef \,h_n\, f\left({\boldsymbol k}, t_n+c_{\ell}\, h_n;\phi_n+\sum_{r=1}^{\ell-1}a_{\ell,r}\,w_r;\mathcal{C}_n\right),
\end{align}
where $a_{\ell,r}$, $b_{\ell}$ and $c_{\ell}$ are the coefficients of the Butcher tableau which characterizes the integrator \citep{wanner1996solving}. 
Using Taylor expansions and un-nesting the scheme, it is possible to prove that a result with a similar structure compared with \eqref{Deltaphi4} is obtained. In this case, the number of stages $s$ appears as exponent in the function $\alpha$, which takes the form  $\alpha({\boldsymbol x}, t;\mathcal{C}_n)=(V_n+
\mathcal{C}_n)^{2s}\,\psi_n$. The function $\beta$, on the other hand, becomes explicitly dependent on $\mathcal{C}_n$, involving a number of terms growing exponentially with $s$.
For this reason, even though the exact result can always be achieved, the computational time needed to minimize the error \eqref{delta_n} is often larger than the time gained with a larger step, especially for higher order schemes ($s>3$).  In the next section, we show how, for practical applications, an exact solution is not necessary to improve the algorithm, and \eqref{L2normres} represents a fast and accurate method.

\section{Numerical examples}\label{secnum}

Here, we consider numerical examples where we apply this method, focusing 
on both the SN and NLS equations solved with the Dormand 
and Prince 5(4) integrator \citep{dormand1980family} in one and two spatial dimensions.
In all 
cases, we set open boundary conditions for the potential, while the initial 
conditions and the value of the physical parameters are chosen to be very 
close to regimes of physical interest, as described in 
\citep{dalfovo1999theory, hasegawa1989optical, kharif2008rogue}. Details 
on the numerical simulations can be found in \ref{num-sim}. 
The gain factor provided by the method depends on the optimal value of 
$\mathcal{C}_n$ compared to the $\mathcal{C}_n=0$ case, which changes from case to case as a function of the boundary conditions for the 
potential and of the profile of the solution. Specifically, since the gain 
factor is evaluated with respect to the $\mathcal{C}_n=0$ case, the more the 
optimal value of $\mathcal{C}_n$ is far from zero, the larger 
the gain factor gets. 

For the one-dimensional NLS, some analytical stationary solutions are known. We then use one of these solutions (see \ref{num-sim}) as initial condition. For all other cases (SN and NLS 2D), no such stationary solutions are known, so we use gaussian initial conditions.

\begin{figure}
	\resizebox{0.47\textwidth}{!}{\includegraphics{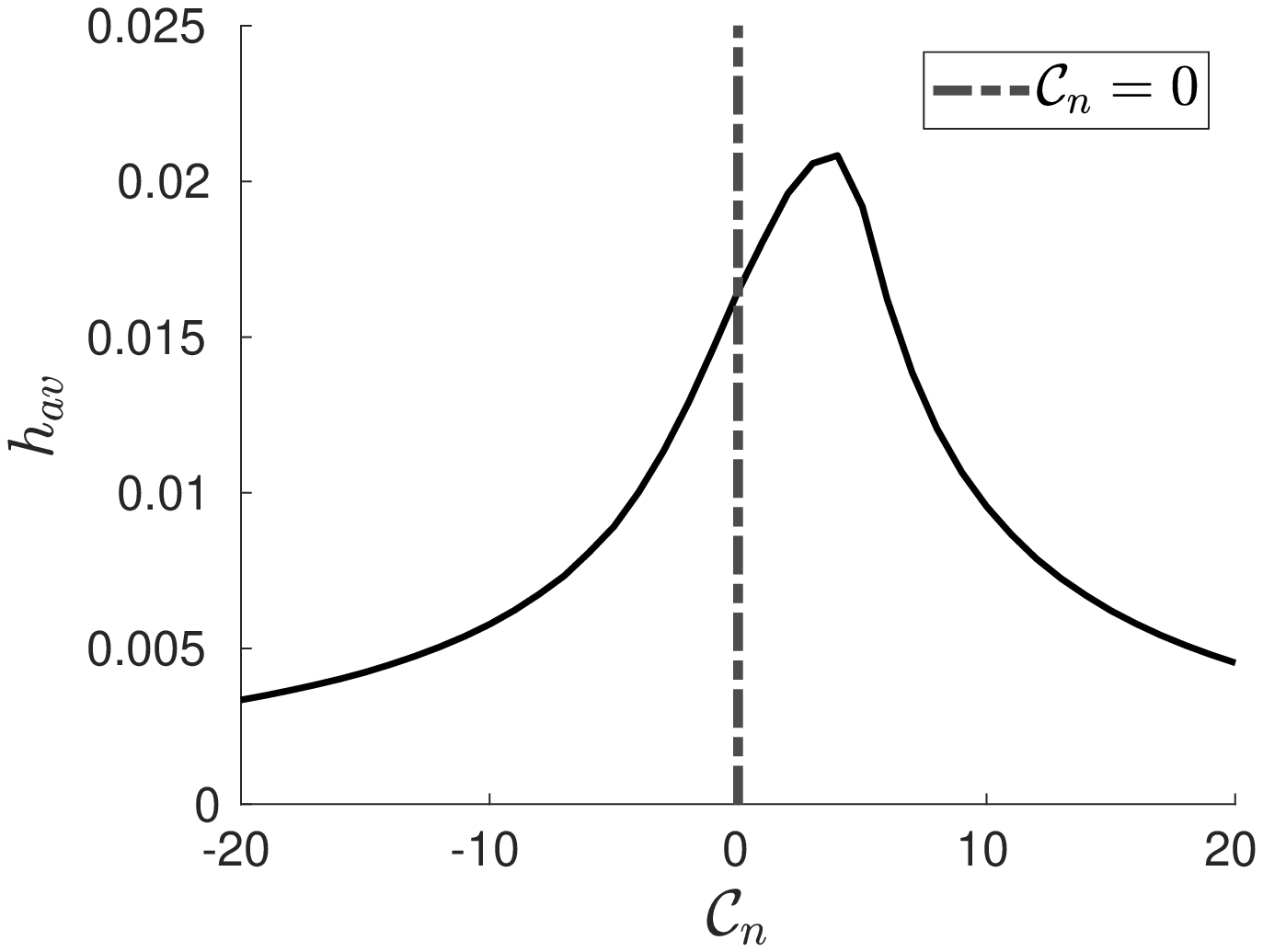}}
	\resizebox{0.47\textwidth}{!}{\includegraphics{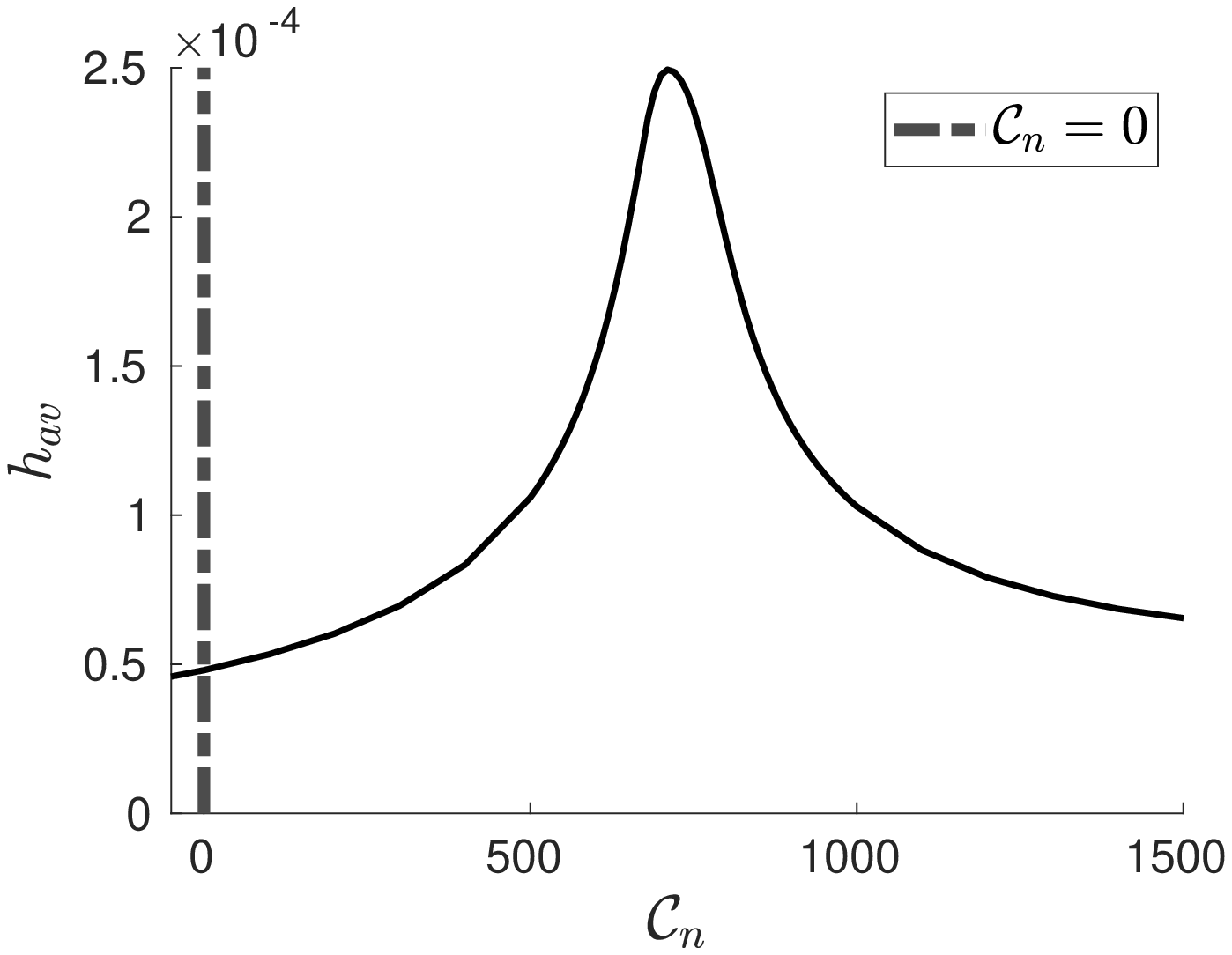}}
	\caption{Average time-step $h_{av}=\textstyle{\frac{1}{N_h}}\sum_{n=1}^{N_h}h_n$ with a constant $\mathcal{C}_n$ for the IF method applied to the one dimensional NLS (left) and SN (right) equations.}
	\label{Coptdtav}
\end{figure}

\begin{figure}
	\resizebox{\textwidth}{!}{\includegraphics{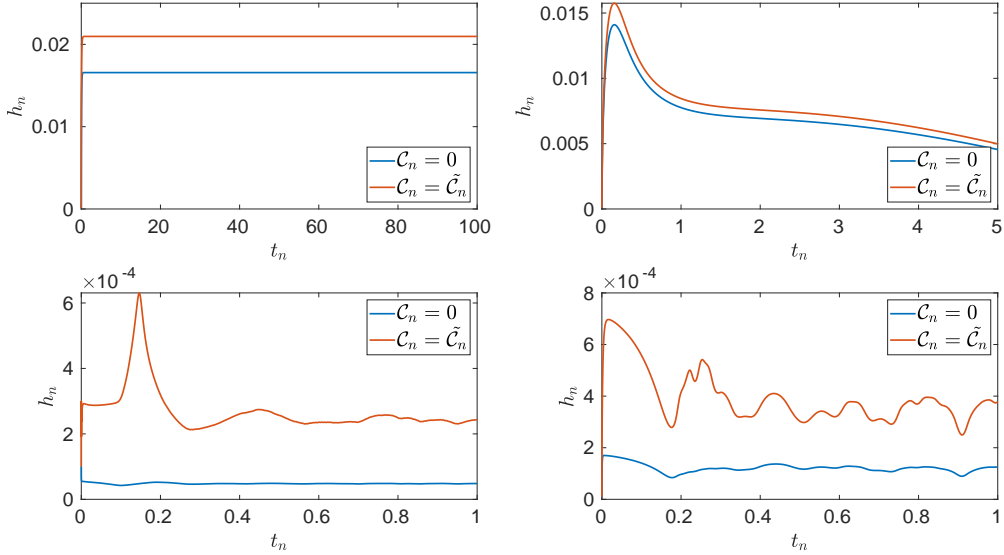}}
	\caption{Comparison between $\mathcal{C}_n=\tilde{\mathcal{C}}_n$ and $\mathcal{C}_n=0$ for time-step $h_n$ as a function of time, for the IF method applied to the $\mathrm{NLS_{1D}}$ (top left),  $\mathrm{NLS_{2D}}$ (top right), $\mathrm{SN_{1D}}$ (bottom left) and $\mathrm{SN_{2D}}$ (bottom down).}
	\label{Coptdtn}
\end{figure}

In Figure \ref{Coptdtav}, we show the average time-step 
$h_{av}=\sum_{n=1}^{N_h}h_n\,/N_h$, for an entire 
simulation with $N_h$ time steps, as a function of 
$\mathcal{C}_n$ for the one-dimensional SN and NLS equations.
These plots are generated taking 
$\mathcal{C}_n$ constant for the entire simulations, in order to better appreciate the strong dependence of the time-step on the choice of the gauge for the potential.
In Figure \ref{Coptdtn}, we report the result of simulations performed 
choosing the near optimal $\mathcal{C}_n=\tilde{\mathcal C}_n$ at each 
time-step. Note that, for the one-dimensional NLS, the solution being 
stationary, $h_n$ and the optimum $\mathcal{C}_n$ do not change in time, 
that is not the case in 2D. We show the time-step $h_n$ as a function of time for the one-dimensional SN and NLS equations, comparing the $\mathcal{C}_n=0$ case with $\mathcal{C}_n=\tilde{\mathcal{C}}_n$. In both cases, the time-step chosen by the algorithm with the optimization of the gauge constant proves to be larger, compared to the $\mathcal{C}_n=0$ case. 
In Table \ref{tabledet}, we show the number of time-loops $N_{\Delta t}$ required to run each simulation and the time $T$ needed to run the simulation (in seconds) for the cases $\mathcal{C}_n=0$ and $\mathcal{C}=\tilde{\mathcal{C}}_n$. 
For NLS, in the one dimensional case we achieve roughly a $30\%$ improvement in terms of speed gain between the $\mathcal{C}=0$ and $\mathcal{C}=\tilde{\mathcal{C}}_n$ cases, while in two dimensions the speed gain is only approximately $10\%$ since, here, the value of $\tilde{\mathcal{C}}_n$ is very close to zero. 
As long as the $SN$ equation is concerned, \eqref{L2normres} proved to reduce remarkably both the number of time loops and the effective time for the total simulation, providing up to a factor $5$ of improvement with respect to the $\mathcal{C}_n=0$ case in 1D and up to a factor $3$ in 2D. 
%
\begin{table}
	\centering
	\begin{tabular}{c|cccccccc} 
		Eq. & $\mathrm{SN_{1D}}$ & $\mathrm{SN_{1D}}$ & $\mathrm{SN_{2D}}$ & $\mathrm{SN_{2D}}$ & $\mathrm{NLS_{1D}}$ & $\mathrm{NLS_{1D}}$ & $\mathrm{NLS_{2D}}$ & $\mathrm{NLS_{2D}}$\\
		$\mathcal{C}$ & $0$ & $\tilde{\mathcal{C}}_n$ & $0$ & $\tilde{\mathcal{C}}_n$ & $0$ & $\tilde{\mathcal{C}}_n$ & $0$ & $\tilde{\mathcal{C}}_n$\\
		$N_{\Delta t}$ & $20819$ & $3871$ & $8382$ & $2682$ & $6047$ & $4781$ & $754$ & $690$\\
		$T(s)$ & $66.7$ & $12.1$ & $12856$ & $4736$ & $18.9$ & $14.5$ & $23769$ & $22843$
	\end{tabular}
	\caption{Comparisons for the SN and the NLS equations, in one and two spatial dimensions, between different values of $\mathcal{C}_n$.}
	\label{tabledet}
\end{table}

\section{Conclusion} \label{sect-conclusions}

Exploiting a gauge condition on the potential, 
we optimized the integrating factor technique applied to the nonlinear 
Schr{\"o}dinger and Schr{\"o}dinger--Newton equations. 
Although the exact values of the piece-wise constant $\mathcal{C}_n$ 
minimizing the error \eqref{delta_n} (therefore maximizing the time-step) 
is in principle always possible to compute (e.g., with a computer algebra 
system), its expression depends on the particular numerical scheme chosen 
and it becomes complicated as the order of the method increases, resulting 
in a high computational cost. However, the near-optimal value obtained from 
the first approach we described, based on the minimization of the $L_2$-norm 
of the nonlinear part of the equation, proved to be an accurate and efficient 
solution in the tested cases. Thus, being computationally extremely cheap 
and independent of the particular numerical scheme employed, this is the 
approach one should choose for most simulations, at least when the computation 
of $\mathcal{N}$ is not very expensive. For Schr\"odinger-like equations 
with hard to compute potentials, most of the computational time is spent 
in the calculation of $\mathcal{N}$. For these very demanding equations, 
the extra cost needed to compute the optimum $\hat{\mathcal{C}}_n$ (instead 
of the near optimum $\tilde{\mathcal{C}}_n$) is negligible in comparison, 
so $\hat{\mathcal{C}}_n$ could be preferable.  

For the cases tested here, we found a speed-up in the computation time up 
to a factor $5$, the speed-up depending on the equation and on the physical 
regime. These examples show that this approach provides significant 
speed improvements, that with minor modifications of the original 
algorithm.
Though we focused on the  nonlinear Schr{\"o}dinger and 
Schr{\"o}dinger--Newton equations, the method principle  
is independent on the particular potential considered, so this approach can 
be extended to other Schr{\"o}dinger-like equations. More generally, the 
idea behind the method presented in this note can be, at least in principle, 
generalized and extended to other equations with similar gauge conditions.

\appendix 

\section{Numerical simulations} \label{num-sim}

For the one-dimensional NLS, we considered the case $g=-1$ (see 
\ref{potNLS}) and we used 
$\psi(x,t=0)=\sqrt{2}\,{\mathrm{sech}}\left(\sqrt{2}\,x\right)$ as initial 
condition. We discretised space with $N=2048$ points, in a computational 
box of length $L=80$. The two-dimensional NLS, which is often employed in optics to model self-focusing 
beams in a medium with a cubic non-linearity 
\cite{zakharov1975nature, konno1979self}, presents a finite time 
(blow-up) singularity \cite{sulem2007nonlinear}. More specifically, whenever 
the initial condition $\psi_0$ satisfies $E_g=\int \ud \boldsymbol{r} \,\psi_0, 
\left(-\frac{1}{2}\nabla^2+\frac{g}{2}\left| \psi_0\right|^2 \right) \psi^*_0 
< 0$, the norm of the solution, or of one of its derivatives,  becomes unbounded 
in finite time. 
For this reason, we stop the simulation at $t_{\mathrm{fin}}=5$, i.e., 
before the singularity occurs. We set $\psi(\boldsymbol{r}, t=0)=\ue^{-r^2/2}/\sqrt{\pi}$ 
as initial condition and we consider the $g=-6$ case, for which the 
corresponding initial energy is $E(g=-6) \approx 0.02$, hence quite close 
to the singular regime; for the spatial discretization we used $L=120$ and 
$N=4096^2$ (squared box with side $L=120$ discretised with $4096
\times4096$ nodes).
For both the one and two dimensional SN equations, we set $g=500$ and 
considered a Gaussian initial condition, $\psi(\boldsymbol{x}, t=0)=\mathcal{N}\ue^{-\left| \boldsymbol{x}\right|^2/2}$ where $\mathcal{N}$ 
is the normalization factor, fixed such that $\int \ud\boldsymbol{x}\left|\psi(\boldsymbol{x}, t=0)\right|^2=1$.
The parameters of the spatial discretization are $L=20$ and $N=2048$ in 
1D, while for the 2D case we set $L=20$ and $N=1024^2$.
To solve the SN and the NLS equation we used the Dormand 
and Prince 5(4) integrator \citep{dormand1980family}.


\begin{thebibliography}{10}
\providecommand{\natexlab}[1]{#1}
\providecommand{\url}[1]{\texttt{#1}}
\expandafter\ifx\csname urlstyle\endcsname\relax
\providecommand{\doi}[1]{doi: #1}\else
\providecommand{\doi}{doi: \begingroup \urlstyle{rm}\Url}\fi

\bibitem{dalfovo1999theory}
Franco Dalfovo, Stefano Giorgini, Lev~P Pitaevskii, and Sandro Stringari.
\newblock Theory of {B}ose--{E}instein condensation in trapped gases.
\newblock {\em Reviews of Modern Physics}, 71(3):463, 1999.

\bibitem{hasegawa1989optical}
Akira Hasegawa.
\newblock Optical solitons in fibers.
\newblock In {\em Optical Solitons in Fibers}, pages 1--74. Springer, 1989.

\bibitem{kharif2008rogue}
Christian Kharif, Efim Pelinovsky, and Alexey Slunyaev.
\newblock {\em Rogue waves in the ocean}.
\newblock Springer Science \& Business Media, 2008.

\bibitem{binney}
James {Binney} and Scott {Tremaine}.
\newblock {\em {Galactic Dynamics: Second Edition}}.
\newblock 2008.

\bibitem{dabby1968thermal}
FW~Dabby and JR~Whinnery.
\newblock Thermal self-focusing of laser beams in lead glasses.
\newblock {\em Applied Physics Letters}, 13(8):284--286, 1968.

\bibitem{rotschild2005solitons}
Carmel Rotschild, Oren Cohen, Ofer Manela, Mordechai Segev, and Tal Carmon.
\newblock Solitons in nonlinear media with an infinite range of nonlocality:
  first observation of coherent elliptic solitons and of vortex-ring solitons.
\newblock {\em Physical review letters}, 95(21):213904, 2005.

\bibitem{rotschild2006long}
Carmel Rotschild, Barak Alfassi, Oren Cohen, and Mordechai Segev.
\newblock Long-range interactions between optical solitons.
\newblock {\em Nature Physics}, 2(11):769--774, 2006.

\bibitem{guzman2006gravitational}
F~Siddhartha Guzman and L~Arturo Urena-Lopez.
\newblock Gravitational cooling of self-gravitating {B}ose condensates.
\newblock {\em The Astrophysical Journal}, 645(2):814, 2006.

\bibitem{hu2000fuzzy}
Wayne Hu, Rennan Barkana, and Andrei Gruzinov.
\newblock Fuzzy cold dark matter: the wave properties of ultralight particles.
\newblock {\em Physical Review Letters}, 85(6):1158, 2000.

\bibitem{paredes2016interference}
Angel Paredes and Humberto Michinel.
\newblock Interference of dark matter solitons and galactic offsets.
\newblock {\em Physics of the Dark Universe}, 12:50--55, 2016.

\bibitem{marsh2019strong}
David~JE Marsh and Jens~C Niemeyer.
\newblock Strong constraints on fuzzy dark matter from ultrafaint dwarf galaxy
  {E}ridanus ii.
\newblock {\em Physical review letters}, 123(5):051103, 2019.

\bibitem{diosi2014gravitation}
Lajos Diosi.
\newblock Gravitation and quantumm-echanical localization of macroobjects.
\newblock {\em arXiv preprint arXiv:1412.0201}, 2014.

\bibitem{penrose1996gravity}
Roger Penrose.
\newblock On gravity's role in quantum state reduction.
\newblock {\em General relativity and gravitation}, 28(5):581--600, 1996.

\bibitem{widrow1993using}
Lawrence~M Widrow and Nick Kaiser.
\newblock Using the {S}chr{\"o}dinger equation to simulate collisionless
  matter.
\newblock {\em The Astrophysical Journal}, 416:L71, 1993.

\bibitem{CanutoEtAl2006-1}
C.~Canuto, M.~Y. Hussaini, A.~Quarteroni, and Th.~A. Zang.
\newblock {\em Spectral Methods. Fundamentals in Single Domains}.
\newblock Scientific Computation. Springer, 2nd edition, 2008.

\bibitem{lawson1967generalized}
J~Douglas Lawson.
\newblock Generalized {R}unge--{K}utta processes for stable systems with large
  {L}ipschitz constants.
\newblock {\em SIAM J. Numer. Anal.}, 4(3):372--380, 1967.

\bibitem{blanes2008splitting}
Sergio Blanes, Fernando Casas, and Ander Murua.
\newblock Splitting and composition methods in the numerical integration of
  differential equations.
\newblock {\em arXiv preprint arXiv:0812.0377}, 2008.

\bibitem{bader2019efficient}
Philipp Bader, Sergio Blanes, Fernando Casas, and Mechthild Thalhammer.
\newblock Efficient time integration methods for {G}ross--{P}itaevskii
  equations with rotation term.
\newblock {\em arXiv preprint arXiv:1910.12097}, 2019.

\bibitem{besse2017high}
Christophe Besse, Genevi{\`e}ve Dujardin, and Ingrid Lacroix-Violet.
\newblock High order exponential integrators for nonlinear {S}chr\"odinger
  equations with application to rotating {B}ose--{E}instein condensates.
\newblock {\em SIAM J. Numer. Anal.}, 55(3):1387--1411, 2017.

\bibitem{alexander1990solving}
Roger Alexander.
\newblock Solving ordinary differential equations i: Nonstiff problems (e.
  hairer, sp norsett, and g. wanner).
\newblock {\em Siam Review}, 32(3):485, 1990.

\bibitem{butcher2016numerical}
John~Charles Butcher.
\newblock {\em Numerical methods for ordinary differential equations}.
\newblock John Wiley \& Sons, 2016.

\bibitem{dormand1980family}
John~R Dormand and Peter~J Prince.
\newblock A family of embedded {R}unge--{K}utta formulae.
\newblock {\em J. Comp. Appl. Math.}, 6(1):19--26, 1980.

\bibitem{suli2003introduction}
Endre S{\"u}li and David~F Mayers.
\newblock {\em An introduction to numerical analysis}.
\newblock Cambridge university press, 2003.

\bibitem{wanner1996solving}
Gerhard Wanner and Ernst Hairer.
\newblock {\em Solving ordinary differential equations II}.
\newblock Springer Berlin Heidelberg, 1996.

\bibitem{gustafsson1991control}
Kjell Gustafsson.
\newblock Control theoretic techniques for stepsize selection in explicit
  {R}unge--{K}utta methods.
\newblock {\em ACM Transactions on Mathematical Software (TOMS)},
  17(4):533--554, 1991.

\bibitem{zakharov1975nature}
VE~Zakharov and VS~Synakh.
\newblock The nature of the self-focusing singularity.
\newblock {\em Zh. Eksp. Teor. Fiz}, 68:940--947, 1975.

\bibitem{konno1979self}
Kimiaki Konno and Hiromitsu Suzuki.
\newblock Self-focussing of laser beam in nonlinear media.
\newblock {\em Physica Scripta}, 20(3-4):382, 1979.

\bibitem{sulem2007nonlinear}
Catherine Sulem and Pierre-Louis Sulem.
\newblock {\em The nonlinear {S}chr{\"o}dinger equation: self-focusing and wave
  collapse}, volume 139.
\newblock Springer Science \& Business Media, 2007.

\end{thebibliography}
\end{document}